\title{ \bf An almost trivial observation about the icosahedron\\
}
\author{ \it J\"{u}rgen Richter-Gebert}

\documentclass[10pt]{article} 

\usepackage[letterpaper,
            bindingoffset=0.2in,
            left=.4in,
            right=1in,
            top=.9in,
            bottom=.9in,
            footskip=.25in]{geometry}
\usepackage{blindtext}
\usepackage{multicol}
\usepackage{lipsum}
\usepackage{cuted}
\usepackage{mwe}
\usepackage{graphicx}
\usepackage{wrapfig}
\usepackage[autostyle=true]{csquotes}

\usepackage{amsthm}
\newtheorem*{theorem*}{Theorem}
\newtheorem*{theoremA*}{Theorem A}
\newtheorem*{theoremB*}{Theorem B}
\newtheorem*{theoremC*}{Theorem C}

\newtheorem{theorem}{Theorem}

\usepackage{amssymb}
\usepackage[font=small,labelfont=bf]{caption} 
\usepackage{color}
\usepackage{mathtools}
\usepackage{amsmath}
\usepackage[dvipsnames]{xcolor}
\usepackage[all]{nowidow}
\usepackage{float}
\usepackage{tikz}

\usepackage{soul} 
\usepackage[normalem]{ulem} 

\usepackage[color links=true, backref=page]{hyperref} 

\hypersetup{
  colorlinks,
  citecolor=blue,
  linkcolor=Red,
  urlcolor=Blue
  }

\begin{document}
\newpage
\maketitle





\def\upto{, \ldots ,}
\def\lines{\square}
\def\points{\bigcirc}
\def\lines{\vee}
\def\points{\wedge}
\parskip=1mm
\def\meet{\wedge}
\def\join{\vee}

\def\pts[#1]{\vee_{#1}}
\def\lns[#1]{\wedge_{#1}}
\def\coa#1{\color[rgb]{.2,.2,1}#1}
\def\cob#1{\color[rgb]{.7,0,.4}#1}
\def\coc#1{\color[rgb]{.8,0,0}#1}
\def\cod#1{\color[rgb]{.8,.6,0}#1}
\def\coe#1{\color[rgb]{0,.6,0}#1}
\def\o#1{\overline{#1}}
\def\u#1{\underline{#1}}

\setlength{\parindent}{.5cm}

%


\parskip=2mm
\begin{multicols}{2}

\section{How rigid is the icosahedron?}
Consider a regular icosahedron with its 12 vertices and 20 faces. As a spherical cell complex requiring only the coplanarity of its triangular faces, it is obviously far from rigid. Moving any of its vertices within a small open neighborhood still creates the same combinatorial type of the convex hull.
In fact by the famous theorem of Steinitz from 1922 every spherical cell complex satisfying mild combinatorial conditions can be realized as polyhedron and its {\it realization space} is topologically an open ball. As a consequence every 3-dimensional polytope can be realized with rational coordinates and any two realizations can be continuously morphed into each other while staying in the same combinatorial type \cite{Z95, RG96}.

The situation changes drastically if one considers 4-dimensional polytopes \cite{RG96}, or if one considers polyhedral complexes of higher genus \cite{BR97}. There are  {\it universality theorems} related to these objects that can encode arbitrary polynomial systems of equations and inequalities into these objects.
As a consequence, disconnected realization spaces, necessarily non-rational coordinates and rigidity forced by coplanarities can occur.

Amazingly, already the vertex set of the icosahedron exhibits all these properties as well, if coplanarities are considered that are related to planes cutting through the interior of the object. Each vertex of the icosahedron has five neighbors, and these neighbors form a pentagon. If we record only the coplanarity information of these twelve pentagons—forgetting all edge lengths, angles, and metric structure—we obtain a purely combinatorial incidence structure on twelve points in space. We will consider realizations of this structure.
The seemingly mild condition of coplanarity of the points in these twelve pentagons turns out to be remarkably restrictive. Up to projective equivalence there exist exactly two realizations of this configuration. Both coincide with the vertex set of the regular icosahedron. The two realizations differ only in how the twelve pentagons are interpreted as faces: up to projective transformations they give rise to the great dodecahedron and the small stellated dodecahedron, two of the classical Kepler–Poinsot star polyhedra. 

In our approach of proving this result we exploit a nice connection to the frequently studied pentagram map that creates one pentagon from another by intersecting the diagonals. Within our configuration we will find a short cycle of such maps forcing $P^2(X)\sim X$ where ``$\sim$'' stands for a homothetic relation (equivalence up to translation and scaling). It turns out that the rigidity of this relation translates into the rigidity of our configuration.

Let me close the introduction with a personal remark: I do not claim that the observation in this paper is new or deep by any means. However, it came as a surprise when I first saw it and to the best of my knowledge I was not able to find any reference in the literature for it. Furthermore, observing that the statement is in tight connection to the {\it Kepler–Poinsot polyhedra} and to {\it spatial incidence theorems} as well as to the {\it pentagram map} immediately raised it in my category of {\it ''things one should know about the icosahedron''.}

\section{The configuration}

 \begin{figure*}[t]
\noindent
\begin{center}
\includegraphics[width=0.97\textwidth]{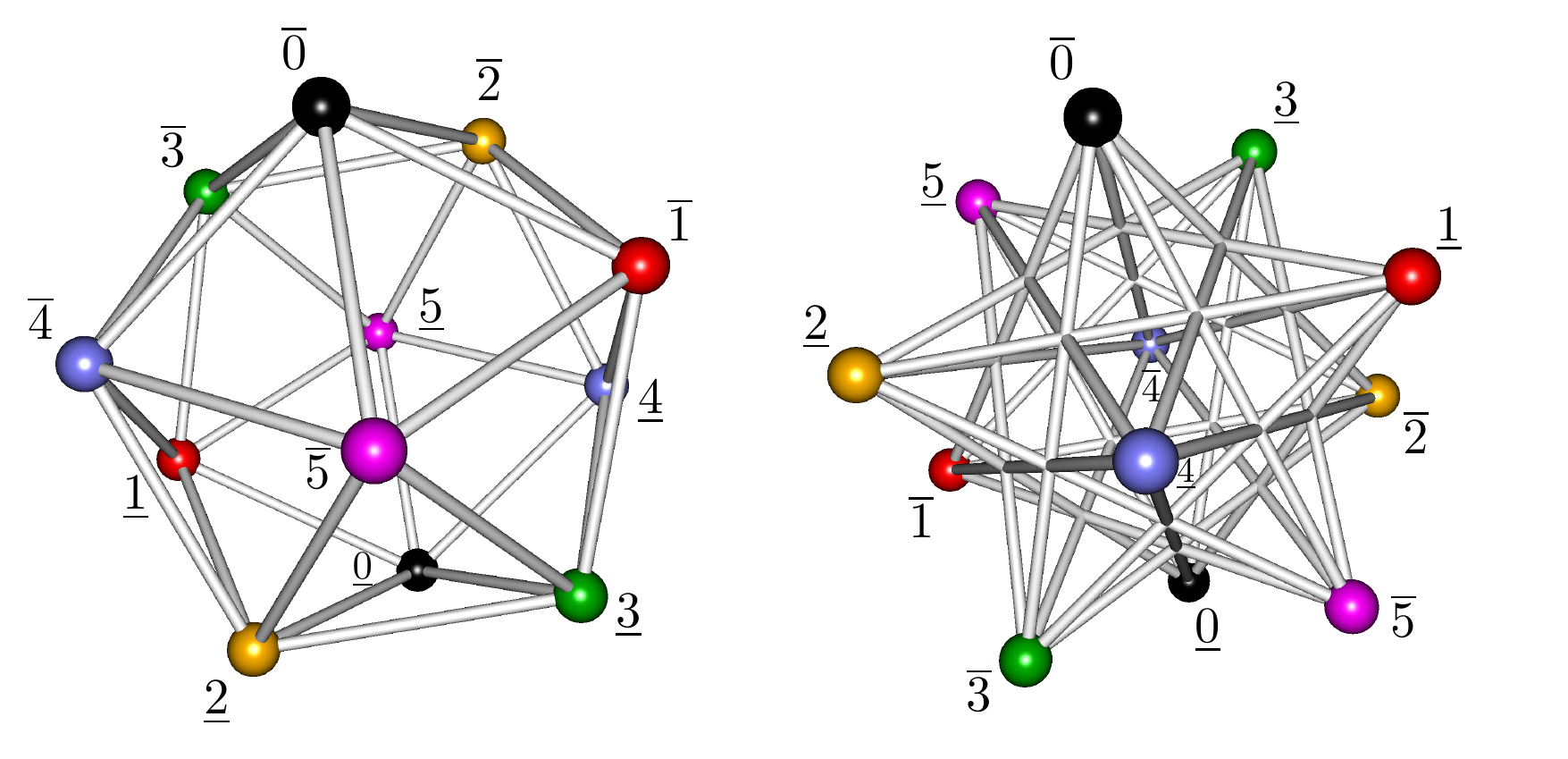}
\end{center}
\captionof{figure}{The edge graph of the great dodecahedron (left) $\mathcal{G}$ is identical to the edge graph of the icosahedron.
The small stellated dodecahedron $\mathcal{G^*}$ (right) also has the same edge graph, however with a geometrically different embedding.
The embedding can be derived from $\mathcal{G}$ by a suitable permutation of the vertices.}  \label{fig:GD_SSD}
\end{figure*}

More formally, the planes relevant for our considerations are the following: Consider the regular icosahedron $\mathcal{I}$ with its set of twelve vertices $V=\{p_0,p_1,\ldots,p_{11}\}$. For each such vertex let the list $N(p)$ be the set of vertices that are connected to that point by an edge of $\mathcal{I}$. Each {\it neighborhood} $N(p)=(q_1,\ldots, q_5)$ of $p$ forms a regular pentagon embedded in $\mathbb{R}^3$. We assume that the points in the list  $N(p)$ are ordered counterclockwise when seen from $p$ (indices modulo~5).

Let $\mathcal{P}=\{ N(p)\  \vert\  p \in V\}$ be the twelve pentagons arising in this way. From this we may derive a purely combinatorial list of coplanarities by setting $D=\{(a,b,c,d,e)\vert (p_a,p_b,p_c,p_d,p_e)\in \mathcal{P}\}$. 
We fix the following concrete indexing of the vertices in which the  twelve  pentagons in $D$ become the 5-tuples.

\[
\begin{array}{r@{\,}c@{\,}lll@{}l}
D&=\big\{&(\o1,\o2,\o3,\o4,\o5),&(\o0,\o5,\u3,\u4,\o2),&(\o0,\o1,\u4,\u5,\o3),\\[2mm]
&&(\o0,\o2,\u5,\u1,\o4),&(\o0,\o3,\u1,\u2,\o5),&(\o0,\o4,\u2,\u3,\o1),\\[2mm]
&&(\u0,\u5,\o2,\o1,\u3),&(\u0,\u1,\o3,\o2,\u4),&(\u0,\o2,\o4,\o3,\u5),\\[2mm]
&&(\u0,\u3,\o5,\o4,\u1),&(\u0,\u4,\o1,\o5,\u2),&(\u1,\u2,\u3,\u4,\u5)&\big\}.\\
\end{array}
\]
Here the points $\o{i}$ and $\u{i}$ are antipodal points and $N(\o0)=(\o1,\o2,\o2,\o3,\o4)$. Compare Figure~\ref{fig:GD_SSD}  on the left for the labelling.
We define the index set: $I=\{ \o0,\ldots, \o5,\u0,\ldots, \u5 \}$. By the usual abuse of notation we will identify the point position in $\mathbb{R}^3$ with the label of the point and speak of the point $i$ rather than $p_i$.

\begin{figure*}[t]
\noindent
\begin{center}
\vspace{-.5cm}
\includegraphics[width=0.15\textwidth]{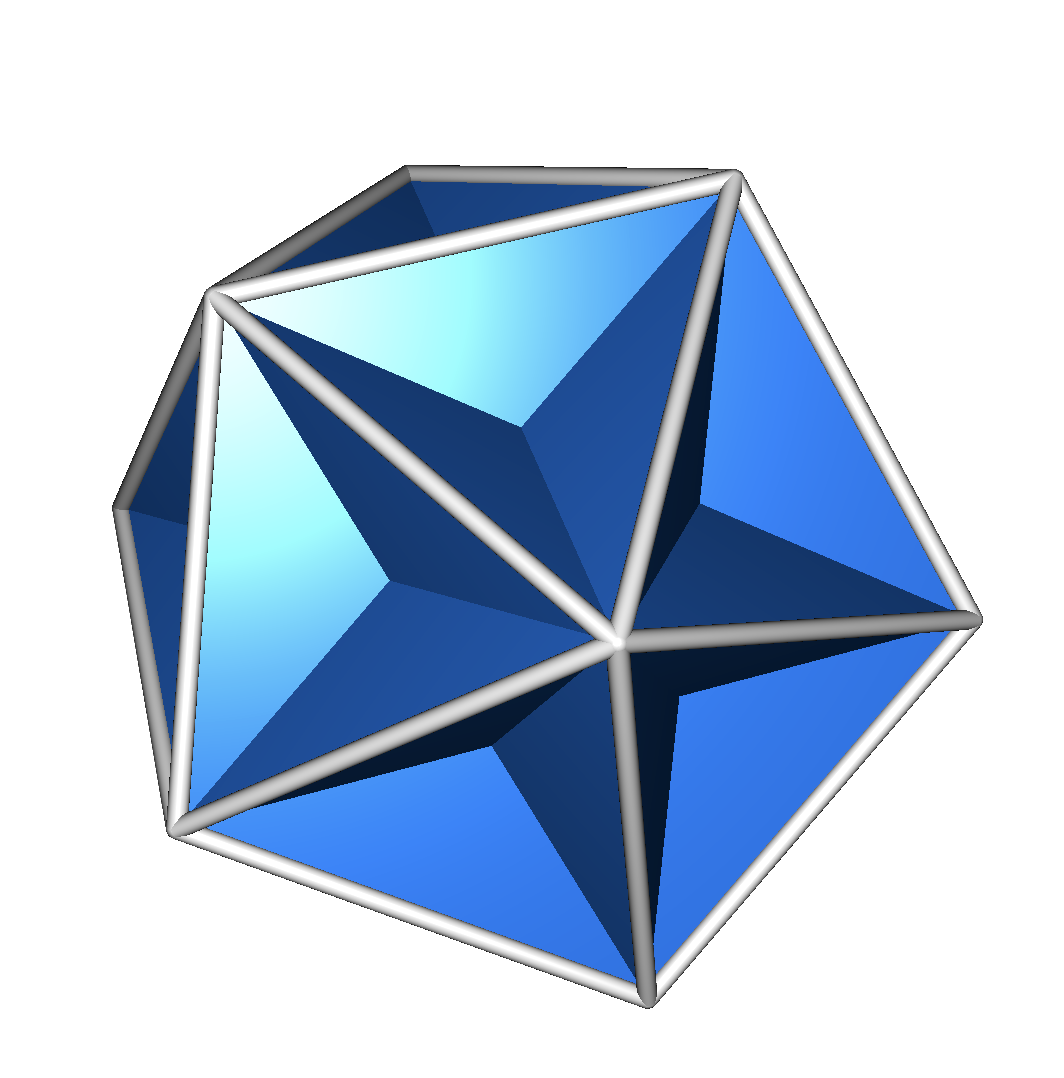}
\raisebox{-4mm}{\includegraphics[width=0.2\textwidth]{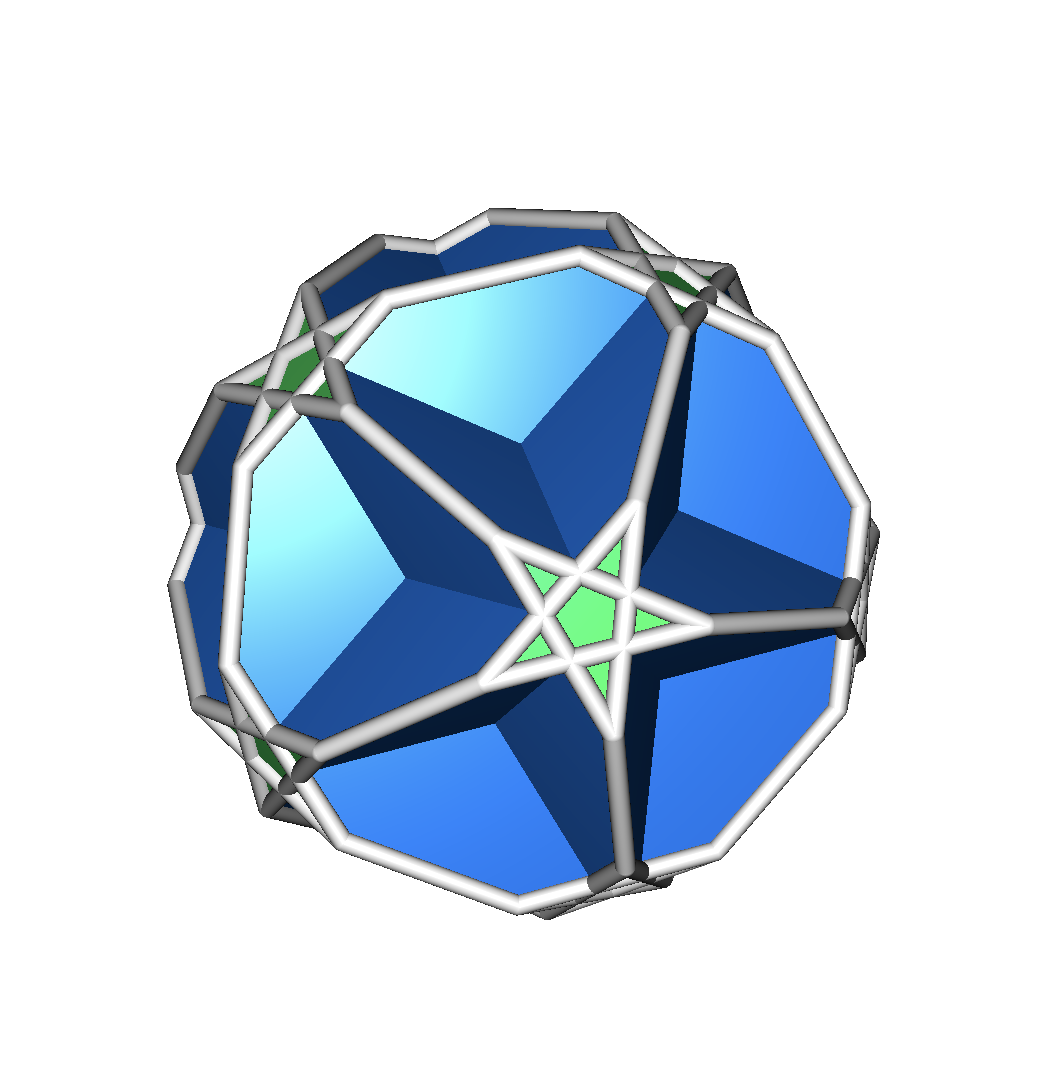}}\hspace{-.4cm}
\raisebox{-8mm}{\includegraphics[width=0.26\textwidth]{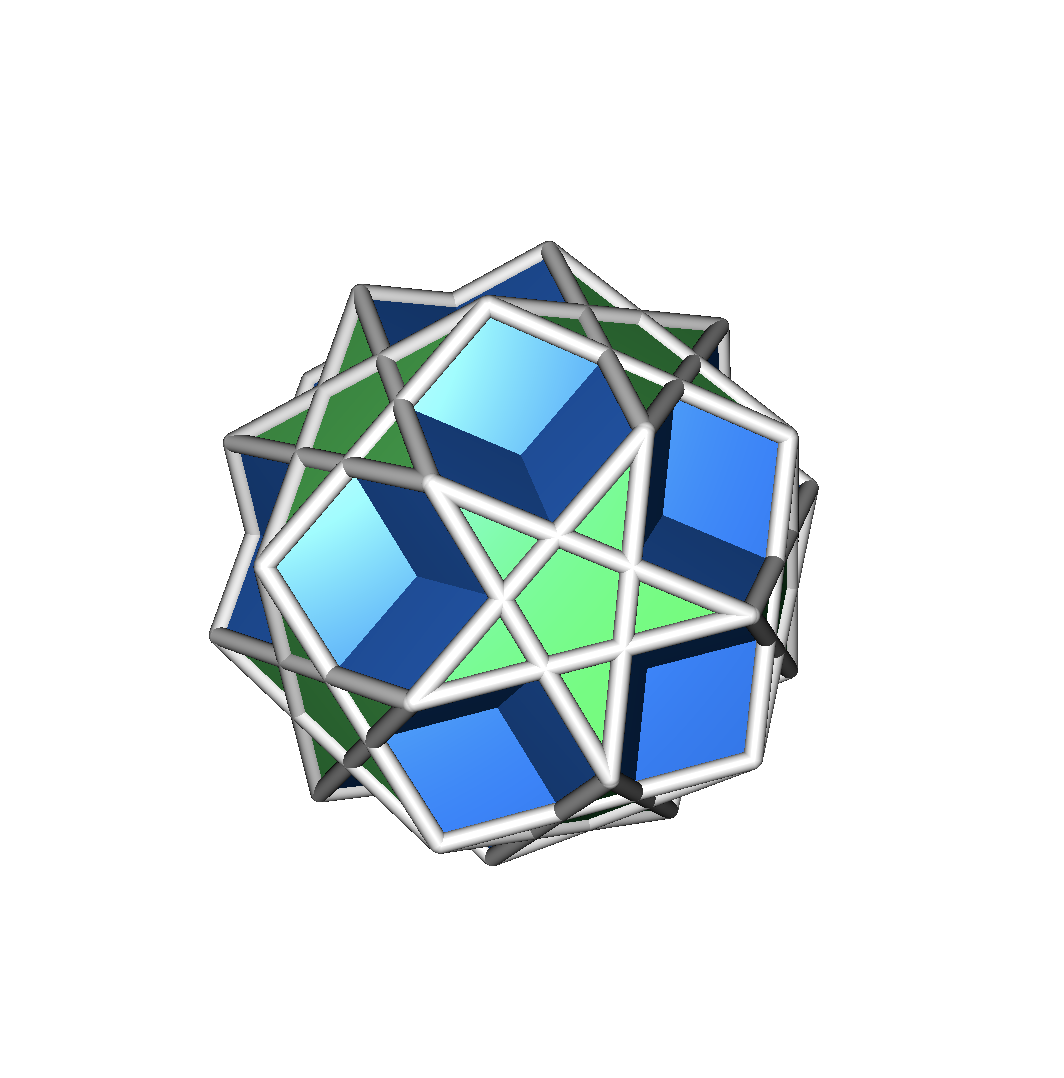}}\hspace{-.4cm}
\raisebox{-4mm}{\includegraphics[width=0.2\textwidth]{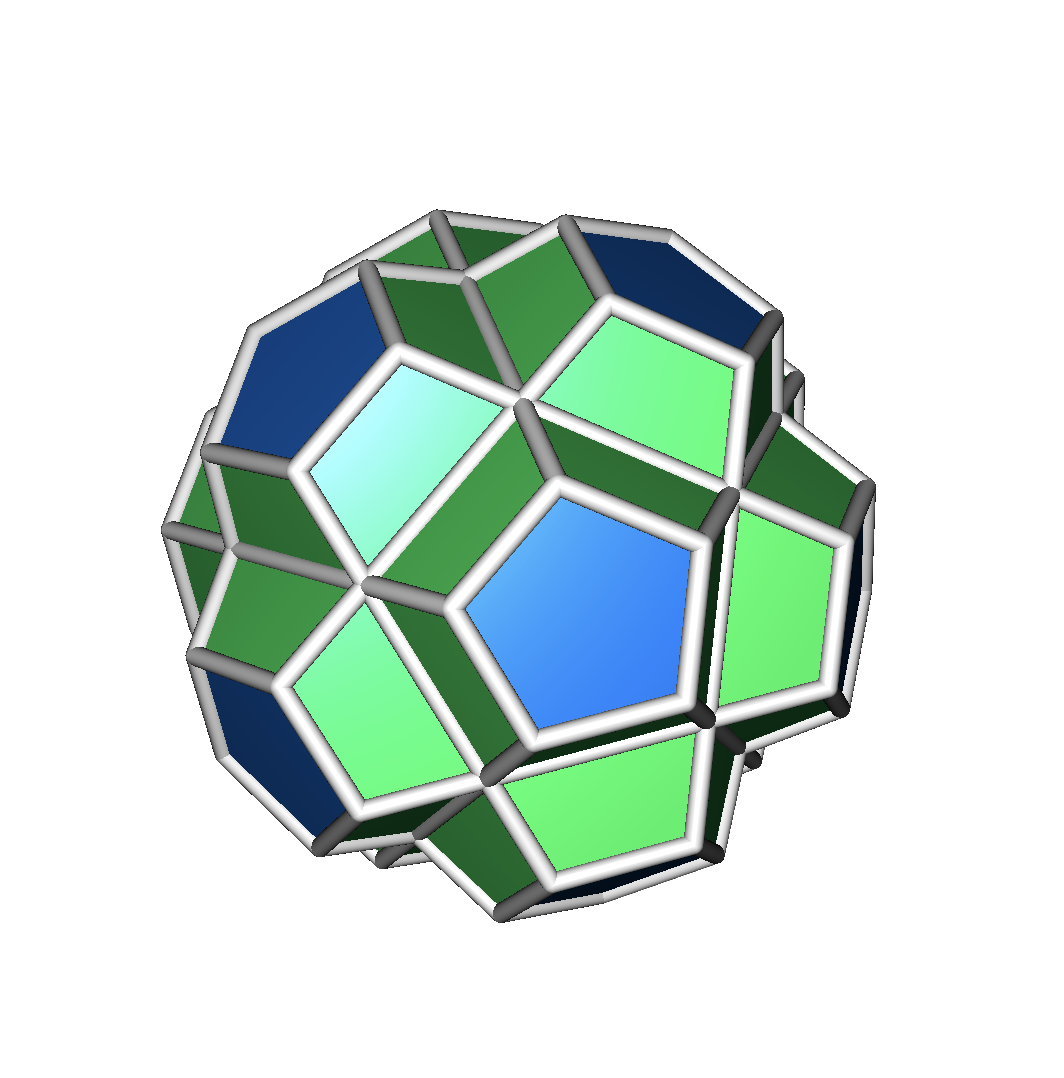}}
\includegraphics[width=0.15\textwidth]{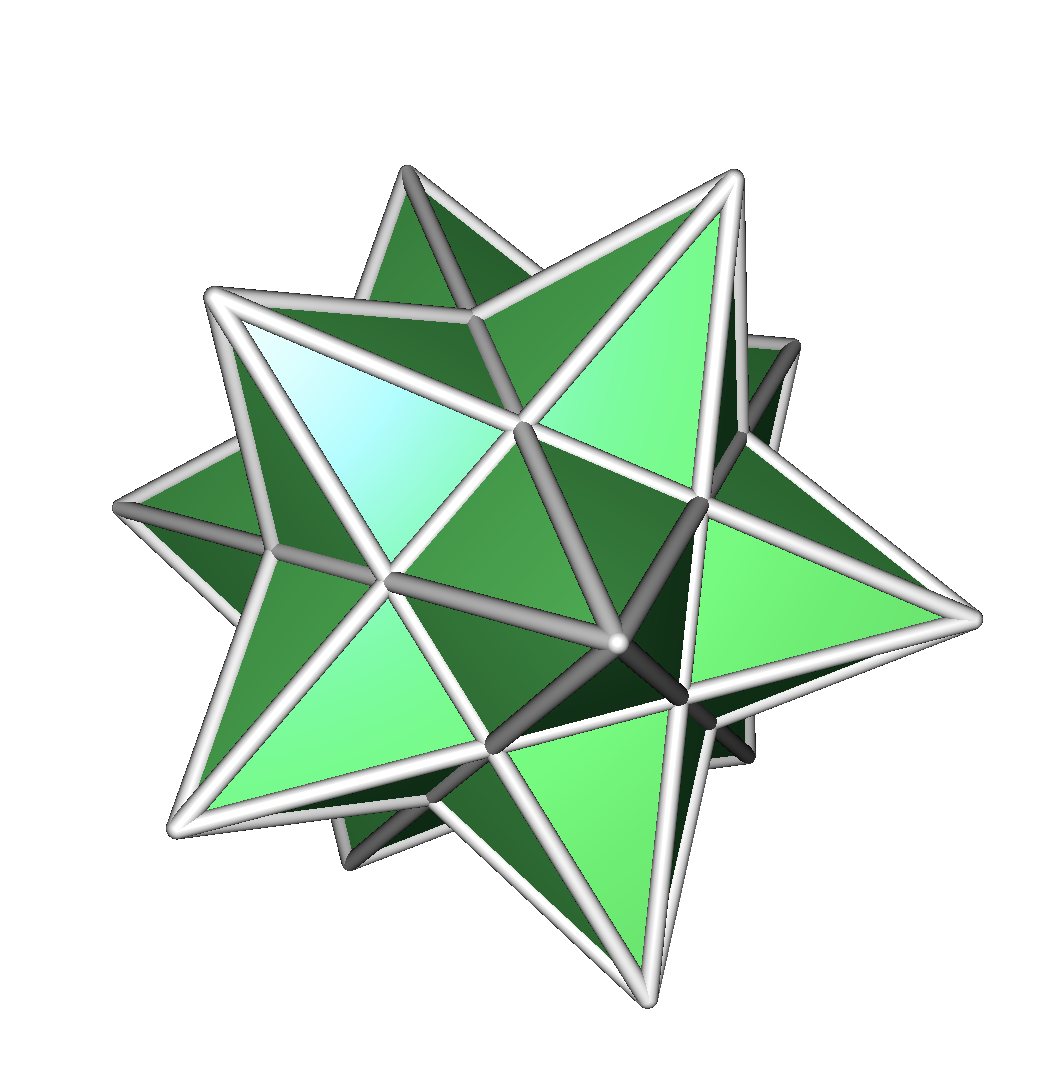}
\end{center}
\vspace{-1cm}
\captionof{figure}{The transition from the great dodecahedron (left) via the dodecadodecahedron (middle) to the small stellated dodecahedron (right)}~\label{fig:DD}
\end{figure*}

\noindent
We will prove:

\begin{theorem}
Consider any set of 12 points $\{ \o0,\ldots, \o5,\u0,\ldots, \u5 \}$ in  $\mathbb{R}^3$ satisfying non-degeneracy $ND$ requirements  specified below. If for each $(a,b,c,d,e)\in D$
the points $a,b,c,d,e$ are coplanar then up to projective equivalence there are only two possibilities for realization.
In both cases after a suitable projective transformation the vertices coincide with the vertex set of $\mathcal{I}$.
Moreover, the second realization is obtained from the first by the following permutation of vertices:
\begin{equation}
\begin{array}{c}
\o0 \rightarrow \o0 ,\ 
\o1 \rightarrow \u1 ,\ 
\o2 \rightarrow \u3 ,\ 
\o3 \rightarrow \u5 ,\ 
\o4 \rightarrow \u2 ,\ 
\o5 \rightarrow \u4 ,\ \\[2mm]
\u0 \rightarrow \u0 ,\ 
\u1 \rightarrow \o1 ,\ 
\u2 \rightarrow \o3 ,\ 
\u3  \rightarrow \o5,\ 
\u4  \rightarrow \o2 ,\ 
\u5  \rightarrow \o4 .\ 
\end{array}
\end{equation}\label{Trans}
\end{theorem}
As non-degeneracy assumptions $ND$ we require that all non-coplanarities of the regular icosahedron also lead to non-coplanar situations in the realization (in other words we are asking for realizations of the {\it matroid} associated to $\mathcal{I}$). In particular, this implies that 
in any realization the planes $A$ and $B$ containing points $\{\o1\upto \o5 \}$ and 
 $\{\u1\upto \u5 \}$, respectively, do not coincide,  points $\o0$ and $\u0$ do not lie in either of these planes and no three points are collinear.
Notice that there are other coplanarities in $\mathcal{I}$ besides the pentagons in $D$. They come from the fact that opposite edges are parallel.

\section{Kepler–Poinsot polyhedra}

Let us consider the twelve pentagons given by the list $D$ as a first class citizen. Each element $(a_1,\upto a_5)\in D$ combinatorially describes a pentagon with edges $a_i,a_{i+1}$ (indizes taken modulo 5). It is easy to check that each edge $(a,b)$ occurs exactly in one pentagon and the opposite edge $(b,a)$ appears in exactly one other pentagon. The twelve pentagons turn out to form one closed surface with 5 pentagons meeting at each vertex. If the pentagons come from the vertex neighbors of the regular icosahedron then the resulting object is known as the great dodecahedron.

The great dodecahedron $\mathcal{G}$ is one of the four Kepler–Poinsot star-solids. These self intersecting polyhedral surfaces share many symmetry properties with Platonic solids (like flag transitivity, constant lengths of Petrie polygons, etc.). The main difference is that they allow for self intersections. The great dodecahedron has 12 faces, 12 vertices and 30 edges and its Euler characteristic $F-E+V=12-30+12=-6$ shows that it is an orientable surface of genus 4.
It is the leftmost object in Figure~\ref{fig:DD}. While its faces are regular pentagons its vertex figures are regular pentagram stars. Figure~\ref{fig:DD} shows a sequence of objects that arise by truncating the vertices of the great dodecahedron. At the end of the sequence (the rightmost image) there is the dual of the great dodecahedron, the so-called small stellated dodecahedron $\mathcal{G^\ast}$. It is a Kepler–Poinsot solid that consists of 12 regular pentagrams. At each vertex 5 of these pentagrams meet forming a pentagonal vertex figure. The pentagrams themselves can be considered as a regular polygon (with self intersections): All edges have same length and all corner angles are identical ($=36^\circ$).

$\mathcal{G}$ and  $\mathcal{G}^*$ are supported by  the same set of vertices, which are the vertices of $\mathcal{I}$. The following fact comes as a little surprise:  While  $G$ and  $G^*$ are geometrically different, they are combinatorially self dual. Creating a list of all pentagons for each of the two solids leads to combinatorially isomorphic structures. The isomorphism between the vertices is given by our permutation in (\ref{Trans}).
 Figure~\ref{fig:GD_SSD} on the left shows the position of edges and vertices after applying the vertex permutation. Notice that the edge graph of 
 $\mathcal{G}$ turns into the edge graph of  $\mathcal{G}^*$. Under the permutation map pairs of antipodal vertices are preserved.

As usual we consider the realization space as the set of all realizations modulo projective transformation.
With this terminology we can reformulate our Theorem 1 and obtain an equivalent statement about Kepler–Poinsot solids:
\begin{theorem}
The realization space of the great dodecahedron contains exactly two points:
The great dodecahedron $\mathcal{G}$ and the small stellated dodecahedron~$\mathcal{G}^*$.
\end{theorem}

\section{A guest appearance of the pentagram map}
The key step is to reinterpret the configuration in a way that reveals an unexpected connection to the pentagram map.
Throughout we have been considering every realization modulo projective transformation. 
This is the standard way to consider realization spaces, since every projective image of a realizations of an incidence configuration turns out to be a realization of the same configuration as well. By this we can normalize our situation and make assumptions that move the situation in a special position. Consider the two planes $A$ and $B$ supporting the two pentagons $(\o1,\ldots,\o5)$ and $(\u1,\ldots,\u5)$, respectively.
Through our non-degeneracy assumptions we know that $A$ and $B$ are not identical. Without loss of generality we may assume that we move their
line of intersection to infinity such that $A$ and $B$ become parallel. Again, by our non-degeneracy assumptions the points $\o0$ and $\u0$ lie outside these two planes.

\begin{figure}[H]
\noindent
\begin{center}
\includegraphics[width=0.4\textwidth]{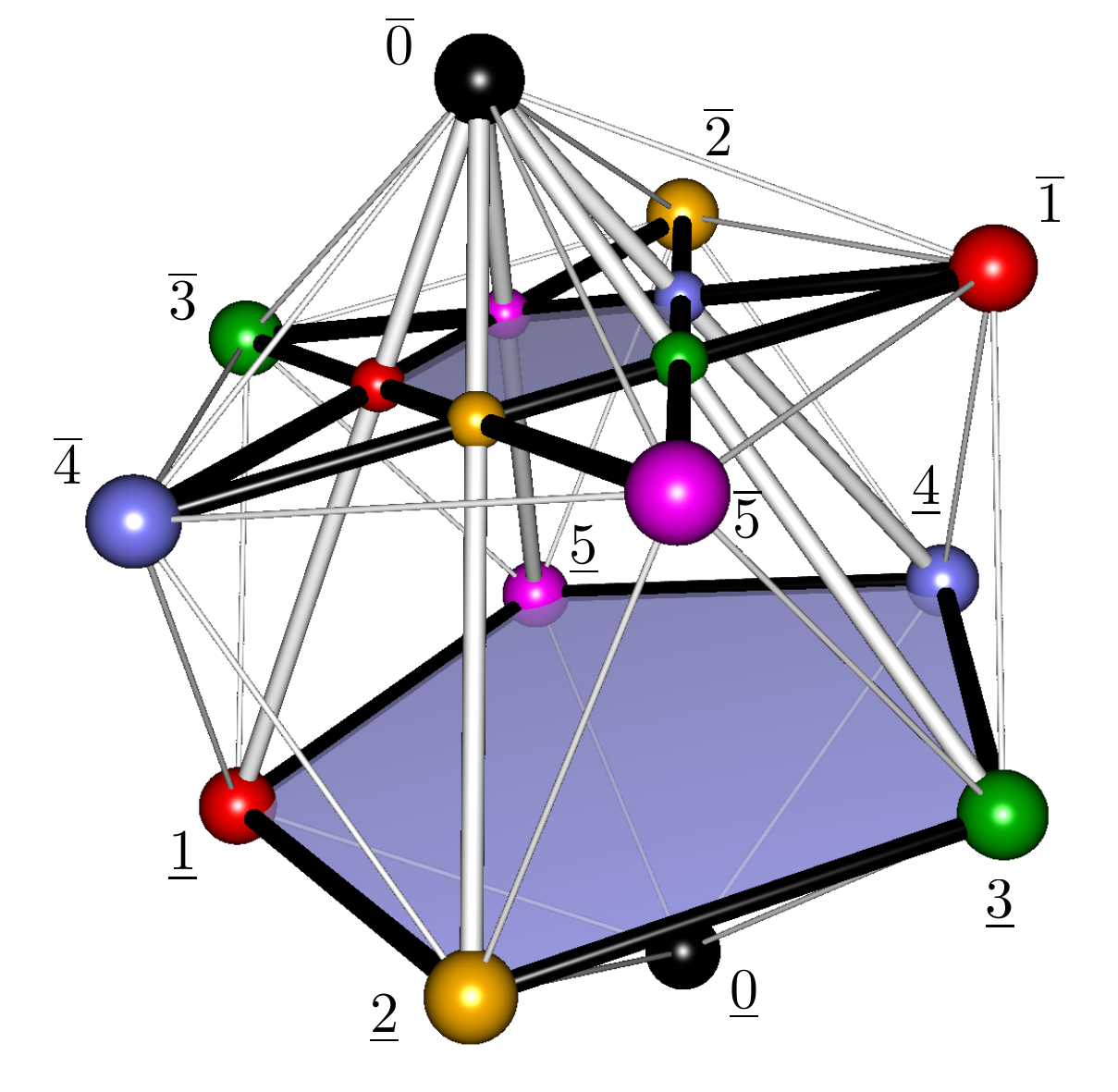}
\end{center}
\captionof{figure}{The projection of $(\o1^*,\ldots,\o5^*)$ to $(\u1,\ldots,\u5)$.}~\label{fig:Project}
\end{figure}

Now consider  Figure~\ref{fig:Project}. It exhibits a sub-configuration of a realization of $D$. The white lines
that meet in point $\o0$ are the intersections of planes $h_{i+1}=\vee(\o0,\o{i},{\o{i+2}})$ and  $h_{i+2}=\vee(\o0,{\o{i+1}},{\o{i+3}})$.
The small colored points are the intersections of these lines with the plane $A$. 
We set 
\[
{\o{i}}^*=\wedge(A,h_{i-1},h_{i+i}).
\]

Hence the points $\o{i}^*$ are the images of the $\u{i}$ when projecting through $\o0$ from $B$ to $A$. 
The points $\o{i}^*$  can be derived 
from $\o{i}$ by an entirely planar construction in the plane $A$: For five points $(1,2,3,4,5)$
we connect points $i$ and  $i+2$ and form their additional five intersections:
\[i^*:=(i+1\vee i+3)\wedge(i-1 \vee i-3)\]
This is the so-called {\it pentagram map} \cite{Sch92}. It associates the points $i^*$ to the points $i$. We denote this map by~$P$.
Hence $P(\o1,\ldots,\o5)=(\o1^*,\ldots,\o5^*)$. Since the planes $A$ and $B$ are parallel  $(\o1^*,\ldots,\o5^*)$ is homothetic (arises by scaling and translation) to $(\u1,\ldots,\u5)$. Reversing the roles of $\o{i}$ and $\u{i}$. We also see that  $P(\u1,\ldots,\u5)$ is homothetic to $(\o1,\ldots,\o5)$.
This allows us to view the entire situation with the plane $A$. Denoting homothety by $\sim$  we just proved 
\begin{theorem}
If $X=(\o1,\ldots,\o5)$ arises from a realization of $D$ then 
\[P^2(X)\sim X.\]
\end{theorem}

It is worth stepping back for a moment and comparing that to other known results about the pentagram map that has been extensively studied in the dynamical systems community and first appears in the work of Alfred Clebsch (1871) \cite{CLE1871}. Already Clebsch proved that for any non-degenerate pentagon $X$ we get a projective equivalence between $X$ and $P(X)$. However, having a homothety between  $X$ and $P^2(X)$ is a significantly stronger geometric requirement that almost completely fixes the possible realizations of $X$.

\begin{figure}[H]
\noindent
\begin{center}
\includegraphics[width=0.48\textwidth]{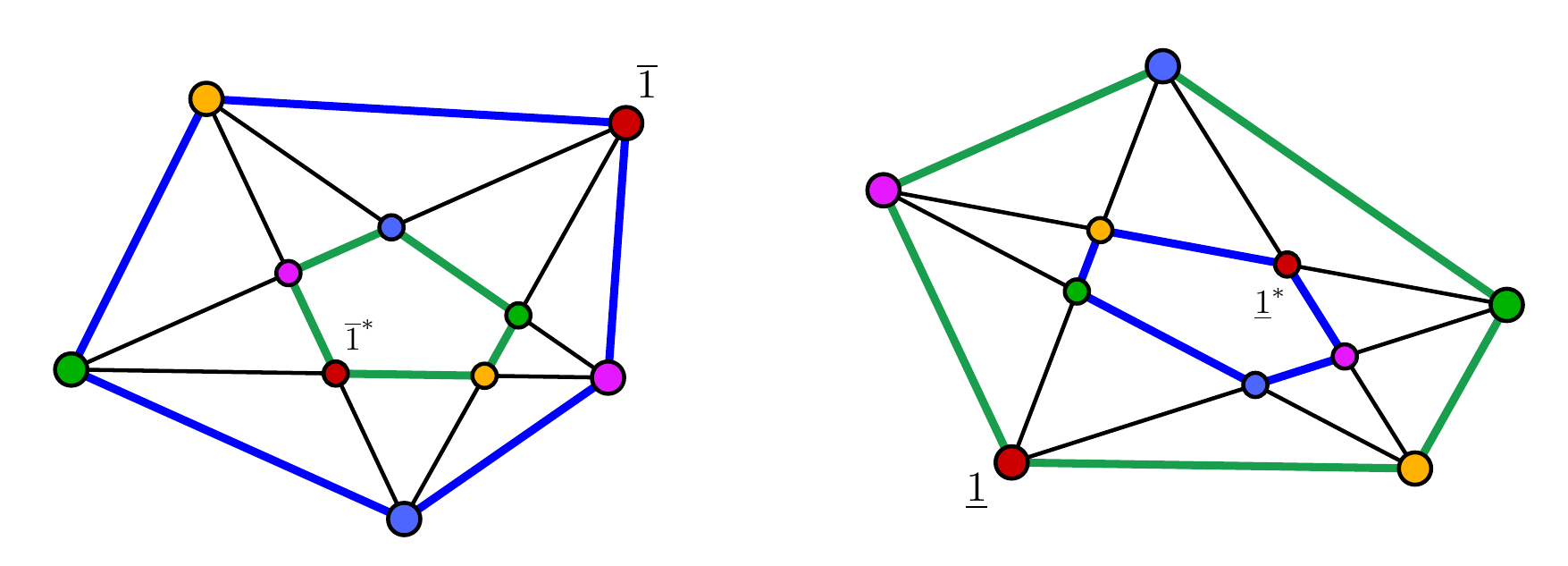}
\end{center}
\captionof{figure}{Pentagram maps between the two pentagons $(\o1,\o2,\o3,\o4,\o5)$ and $(\u1,\u2,\u3,\u4,\u5)$.}  \label{fig:Penta1}
\end{figure}

\noindent
We get 
\begin{theorem}
If a pentagon $X$ satisfies  $P^2(X)\sim X$, then $X$ is either the affine image of a regular pentagon or of a regular star pentagram.
\end{theorem}
After all we said, proving this statement immediately implies Theorem 1 and Theorem 2 as corollaries.

\begin{figure}[H]
\noindent
\begin{center}
\includegraphics[width=0.4\textwidth]{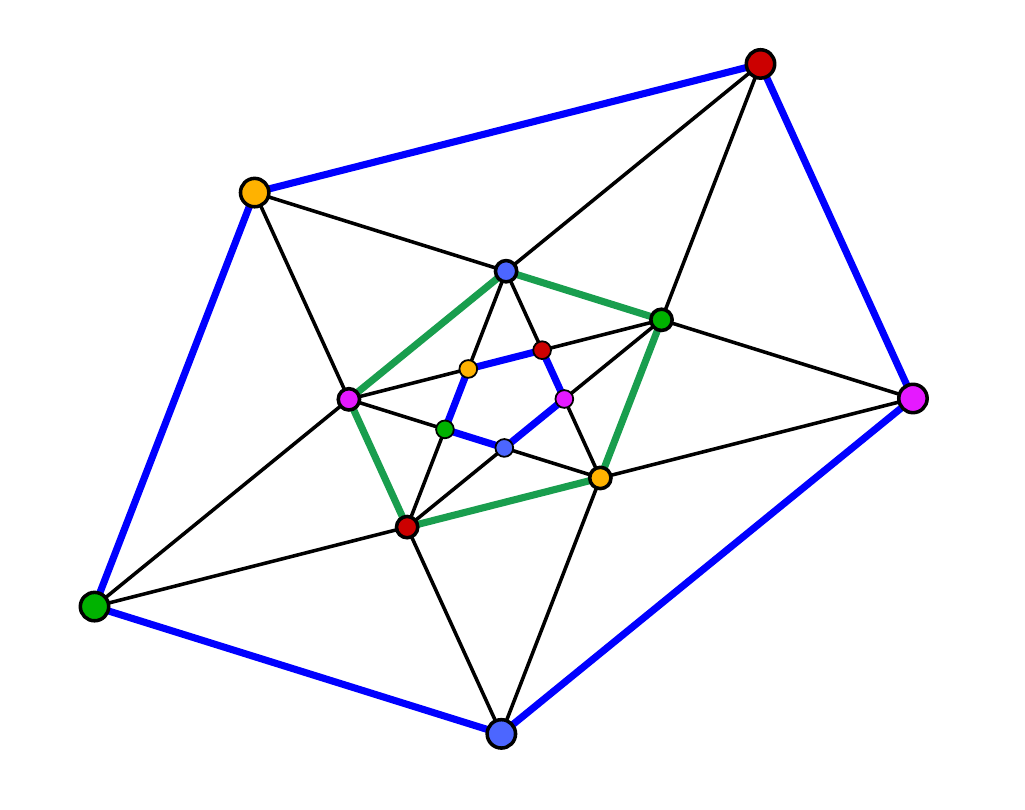}
\end{center}
\captionof{figure}{A geometric situation in which the doubly iterated pentagram map produces a homothetic copy of the original pentagon.}  \label{fig:Penta2}
\end{figure}

Figure~\ref{fig:Penta1} illustrates the situation within the geometry of the planes $A$ and $B$. The left image shows how we start with a blue pentagon $X=(\o1,\o2,\o3,\o4,\o5) $ and construct its pentagram map image $P(X)$. This is homothetic to $Y=(\u1,\u2,\u3,\u4,\u5)$ from which we construct its pentagram map image $P(Y)$ in the right picture. The geometric requirement that this should be homothetic to $X$ is not satisfied in that specific image, hence the pentagons do {\it not} arise from a realization of the great dodecahedron.

Figure~\ref{fig:Penta2} in turn shows a situation where the requirement is satisfied. We put the entire construction into one image. We start with $X$ (blue) construct its pentagram map image $P(X)$ (green) and from that an iterated pentagram map $P^2(X)$ (blue again). In this image the inner blue 
pentagon is a scaled copy of the outer blue pentagon (notice the parallelism of corresponding sides). From there we will now prove that the vertices of $X$ necessarily must, up to affine transformations, be the points of a regular pentagon or regular pentagram star.

\section{The nitty gritty calculations}
Now let's roll up the sleeves and perform the calculations that prove
the rigidity of the configuration. In essence everything boils down to transferring the problem
to a polynomial equation. 
Since we expect exactly two solutions, it is natural to expect a quadratic equation.
Although it is possible to do the calculations in a completely coordinate free invariant theoretic setting we will perform them in a specific coordinate system, to make the considerations as elementary as possible. 

\begin{figure*}[t]
\noindent
\begin{center}
    \begin{tikzpicture}
        \draw (0, 0) node[inner sep=0] {\includegraphics[width=7cm]{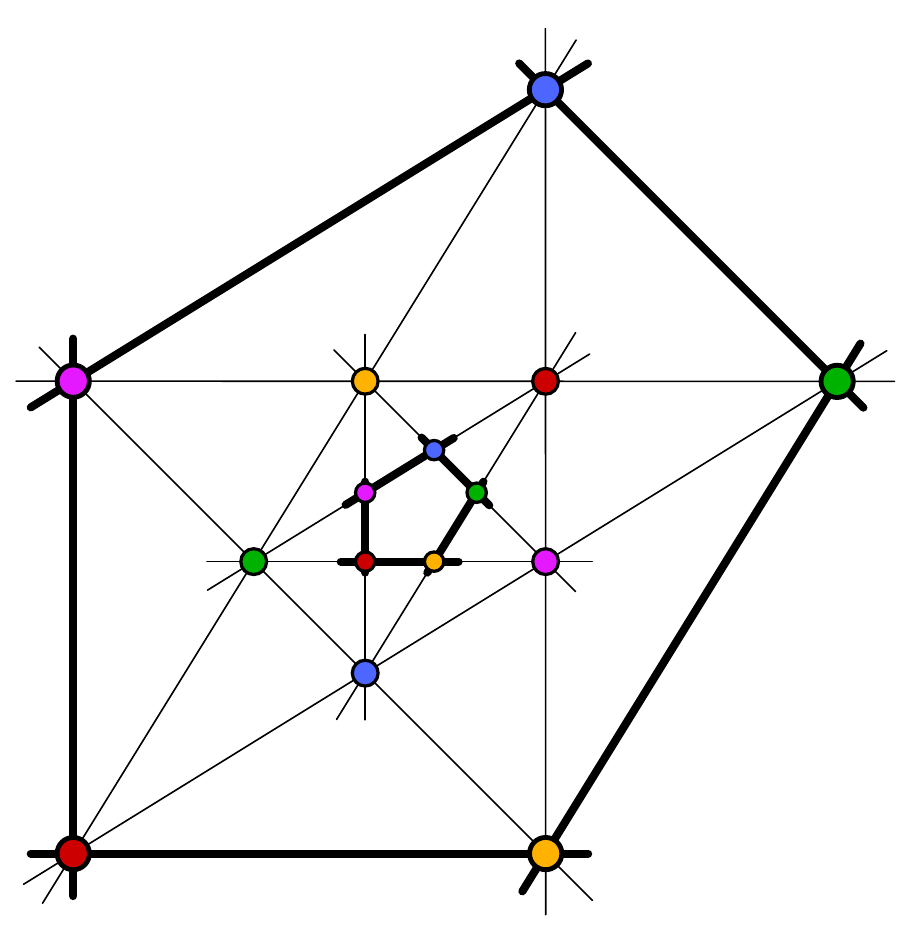}};
        \draw (-9, 0.15) node[inner sep=0] {\includegraphics[width=9cm]{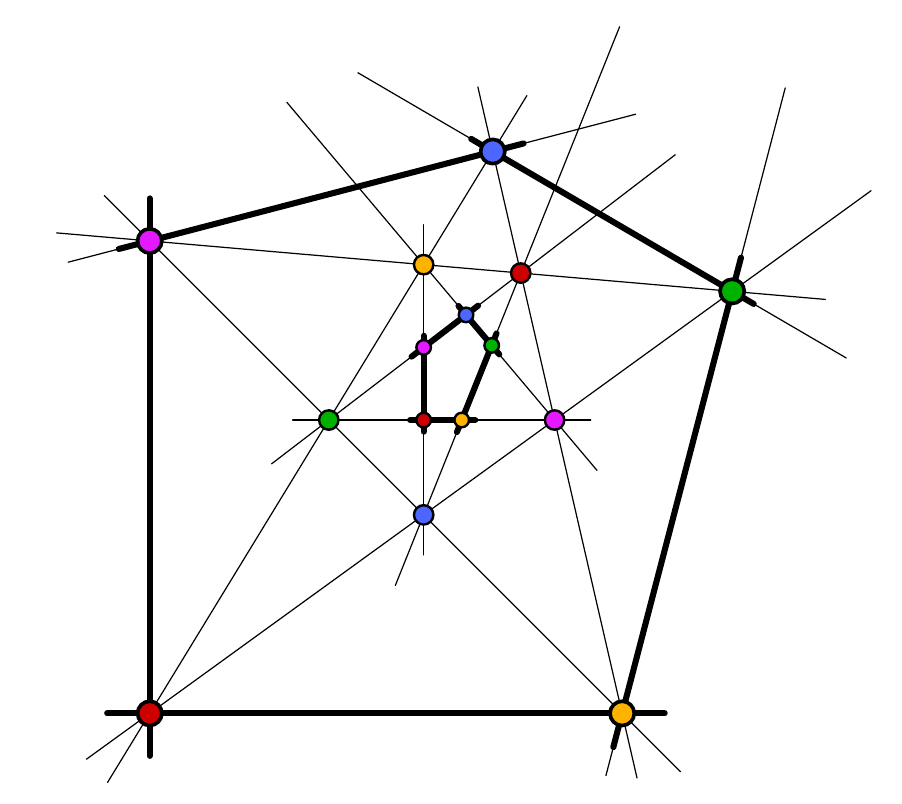}};

         \node at (-11.5,-3.2) {$(0,0)$};
         \node at (-7.9,-3.2) {$(1,0)$};
          \node at (-11.5,2.2) {$(0,1)$};

         \node at (-12.3,-2.7) {$p_1$};
         \node at (-6.9,-2.7) {$p_2$};
         \node at (-6,.7) {$p_3$};
         \node at (-8.9,3.1) {$p_4$};
         \node at (-12.3,1.4) {$p_5$};

         \node at (-8.0,1.2) {$q_1$};
         \node at (-9.6,1.4) {$q_2$};
         \node at (-10.2,-.4) {$q_3$};
         \node at (-8.9,-1.) {$q_4$};
         \node at (-8.1,-0.4) {$q_5$};
         
         \node at (-9.45,-0.3) {$r_1$};

         \node at (-3.2,-2.7) {$p_1$};
         \node at (1.3,-2.7) {$p_2$};
         \node at (3.1,.2) {$p_3$};
         \node at (.3,3.) {$p_4$};
         \node at (-3.2,1.1) {$p_5$};

         \node at (1.0,0.4) {$q_1$};
         \node at (-1.1,0.4) {$q_2$};
         \node at (-1.5,-1.1) {$q_3$};
         \node at (-.3,-1.65) {$q_4$};
         \node at (0.6,-1.1) {$q_5$};
         
          \node at (-0.85,-1.0) {$r_1$};
          
                   \node at (-1.75,-3.35) {$\underbrace{\kern 2.2cm}_{0.618\ldots}$};

    \end{tikzpicture}
\end{center}
\captionof{figure}{A construction sequence for the iterated pentagram map.}  \label{fig:Constr}
\end{figure*}
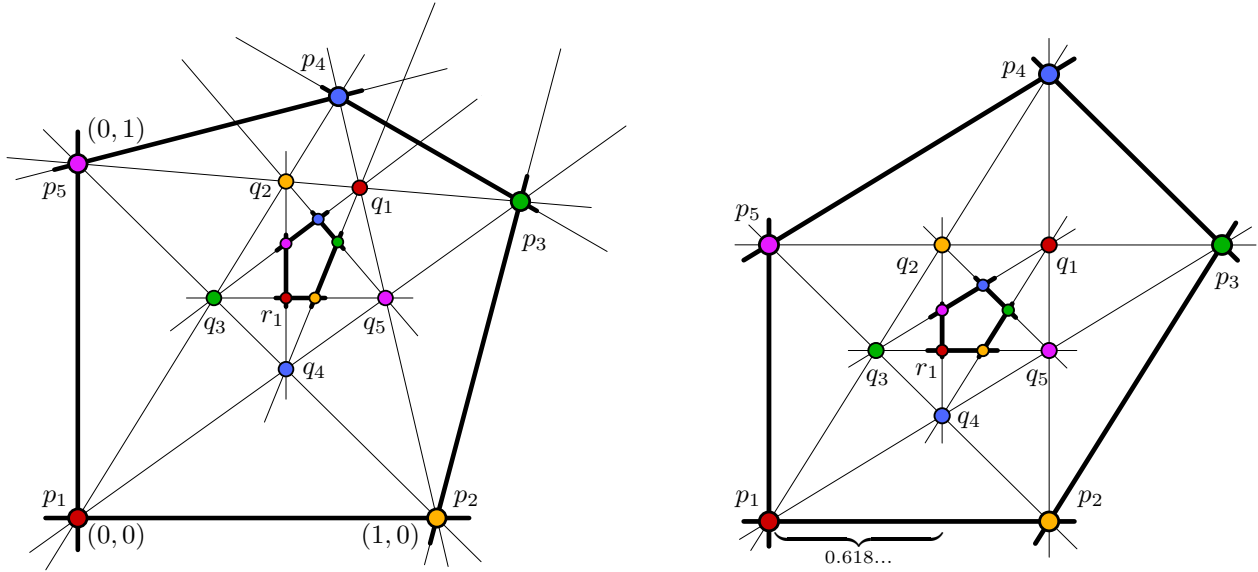

We do the calculation by assigning a special affine basis $(0,0), (1,0),(0,1)$ to the points $1,2,5$, respectively and construct positions of the remaining points by use of {\it join}, {\it meet} and {\it parallel}
operations. All calculations are very elementary. We perform the calculations in homogeneous coordinates. Join ($\vee$) and meet ($\wedge$) can be carried out by vector products. Parallels can be constructed by using intersections with the line at infinity $l_\infty$. We refer to Figure~\ref{fig:Constr} for the labelling, there only the points which are relevant for the calculation are labelled. Besides the affine basis $p_1=(0,0,1),p_2=(1,0,1),p_5=(0,1,1)$, we also introduce the points at infinity in $x$ and $y$ direction
$\infty_x=l_\infty\wedge(p_1\vee p_2)=(1,0,0)$ and
$\infty_y=l_\infty\wedge(p_1\vee p_5)=(0,1,0)$.
We also introduce the point $r_1=(x,y,1)$ as free point from which we calculate positions of all remaining points.

\noindent
 We first calculate:
\[
\begin{array}{rcll}
q_3&=&(p_2\vee p_5)\wedge(r_1\vee \infty_x)=(x, 1 - x, 1),\\[1mm]
q_4&=&(p_2\vee p_5)\wedge(r_1\vee \infty_y)=(1-y, y, 1),\\[1mm]
q_2&=&(p_1\vee q_3)\wedge(r_1\vee \infty_x)=(x - x y, x y, 1 - y),\\[1mm]
q_5&=&(p_1\vee q_4)\wedge(r_1\vee \infty_y)=(x y, y - x y, 1 - x),\\[1mm]
p_4&=&(p_1\vee q_3)\wedge(p_2\vee q_5)\\[1mm]&=&((x-1) (y-1), (1-x)y, 2 - 2 x - y)y,\\[1mm]
p_3&=&(p_1\vee q_4)\wedge(p_5\vee q_2)\\[1mm]&=&((1-y)x, (x-1) (y-1), 2 - x - 2 y)x,\\
\end{array}
\]

Let us analyse how the non-degeneracy assumptions interact with these equations. The non-degeneracy assumptions imply that in  the quintuple of points $(p_1\upto p_5)$  no three points are collinear. 

\noindent
Considering  
\[\det(p_3,p_4,p_5)=x y (1 - 2 x) ( y-1) ( x + y-1)\]
we see that since these three points are not allowed to become collinear none of the factors in that expression may vanish. By considering $\det(p_2,p_3,p_4)\neq 0$ we similarly conclude that in addition
$1 - 2 y\neq 0$ and $x-1\neq 0$.



We now will take care of the first parallelism that is not yet covered by our construction. Line $l_1=q_2\vee q_5$ should be parallel to line $l_2= p_3\vee p_4$ which can be forced by setting $\det(l_1,l_2,l_\infty)=0$.
This expression factorizes to
\[\det(l_1,l_2,l_\infty)=2 {(x - y)}  x y(x-1) ( y-1) ( x + y-1)^2\]
 The only factor that is not among our non-degeneracy conditions and by this does  not force a complete collapse of the configuration is $x-y$.
 This implies that under our non-degeneracy assumptions point $r_1$ lies on the $x=y$ line and the entire configuration becomes symmetric
 w.r.t this axes. With the setting $r_1=(x,x,1)$ we now calculate:
  \[
\begin{array}{rcll}
q_1&=&(p_5\vee q_2)\wedge(p_2\vee q_5)\\[1mm]
&\sim&((x-1 ) x, (x-1 ) x,  2 x^2-1),\\[1mm]
\end{array}
\]
Here we pulled out a non-zero factor of $2x-1$. We now require that the lines $g_1=q_3\vee q_1$ and  $g_2=p_5\vee p_4$ are parallel
by checking $\det(g_1,g_2,l_\infty)=0$. We get:
\[\det(g_1,g_2,l_\infty)=-2 (-1 + x) (-1 + 2 x) (-1 + x + x^2).\]
The only relevant factor is $-1 + x + x^2$ setting this to~$0$ finally leads to:
\[
x={-1\pm\sqrt{5}\over 2}.
\]
The final missing parallelism between 
$q_4\vee q_1$ and  $p_1\vee p_2$
is automatically satisfied by symmetry. In Figure~\ref{fig:Constr} on the right the construction sequence for a generic point $r_1$ is shown, making it visible that generically the three parallelism requirements are not satisfied. The image on the right shows one of the two correct positions.

\medskip
In essence this calculation shows that there are (up to projective equivalence)
only two realizations for $P^2(X)\sim X$. Which in turn shows us that the only realizations we may get for our three dimensional configuration $D$ are the great dodecahedron and the small stellated dodecahedron.


\section{Incidence Theorems}
We close our investigations on the icosahedron with a couple of observations that are a consequence of the rigidity of the configuration $D$.
First of all we perform a degree of freedom count. Our configuration contains 12 points and a projective transformation in $\mathbb{R}^3$ makes it possible to move 5 points to a specified position.
Hence, the a generic configuration has $(12-5)\cdot 3=21$ Degrees of freedom. On the other hand we have $12$ coplanar quintuples. Since a plane is determined by 3 points we have $(5-3)\cdot 12=24$ non-trivial incidences.

We proved that the realization space of the incidence structure is 0-dimensional (just two points). Each degree of freedom can be used to achieve one incidence. So the degree of freedom count gives $21-24=-3$ degrees of freedom which means that three of the incidences are a consequence of the other 21. We have not fully analysed the exact interdependences of the incidences yet.

Another point is remarkable as well. Besides the incidences associated to the pentagons there are also 15  other non-trivial coincidences that are present in any realization of $D$. This can be immediately seen since we know the exact realizations come from the vertex set of the regular icosahedron: In every realization antipodal edges are coplanar. This causes the coplanarity of certain quadruples of points.
It also implies that all connections of antipodal point pairs meet in one single point (the center of the icosahedron).

This statement cannot be reversed. One can linearly morph from the realization of
$\mathcal{G}$ to the realization of $\mathcal{G}^*$ by moving the vertices along a linear path.
This transition can be organised in way that all pairs of opposite edges remain parallel. However, during  the transition all coplanarities of the 15 pentagons get distorted.

\bibliographystyle{plain} 
{\small
\bibliography{refsd}
}
\medskip

\noindent

\noindent
{\footnotesize
\noindent
{\sc Department of Mathematics, Technical University of Munich, Germany}\\
Email address: {\tt richter@tum.de}
}

\end{multicols}

%
%
%
%
%
%
%

\end{document}